\documentclass[12pt]{article}
\usepackage{amsfonts,latexsym,amssymb}
 \newcommand{\beqn}{\begin{eqnarray}}
 \newcommand{\eeqn}{\end{eqnarray}}
 \newcommand{\be}{\begin{equation}}
 \newcommand{\ee}{\end{equation}}
 \newcommand{\ba}{\begin{array}}
 \newcommand{\ea}{\end{array}}
 \newcommand{\pa}{\partial}
 
 \newcommand{\ci}{\cite}
 \newcommand{\la}{\label}

\newcommand{\Om}{\Omega}
\newcommand{\lb}{\lambda}
\newcommand{\ga}{\gamma}

\newcommand{\na}{\nabla}
\newcommand{\om}{\omega}
\newcommand{\bt}{\beta}
\newcommand{\al}{\alpha}

 \newcommand{\De}{\Delta}

\newcommand{\C}{\mathcal{C}}
\newcommand{\rb}{\mathbb{R}}
\newcommand{\wq}{\widetilde{Q}}
\newcommand{\wcl}{\widetilde{\cl}}
\newcommand{\cl}{\mathcal{P}}
\newcommand{\A}{\mathcal{A}}

\newcommand{\wal}{\widetilde{\alpha}}
\newcommand{\M}{\mathcal{M}}

\def\R{{\rm I\kern-.1567em R}}
\def\M{{\rm I\kern-.1567em M}}

\def\div {{\rm div}}

\def\spt{{\rm spt}}

\def\ess{{\rm ess}}

 \newtheorem{theorem}{Theorem}[section]
 
 \newtheorem{definition}[theorem]{Definition}
 
 \newtheorem{lemma}[theorem]{Lemma}

 \newtheorem{cor}[theorem]{Corollary}
 \newtheorem{pro}[theorem]{Proposition}

\begin{document}

\begin{center} {\Large\bf  A sufficient condition of regularity
for axially symmetric solutions to the Navier-Stokes equations }\\
  \vspace{1cm}
 {\large
  G. Seregin and  W. Zajaczkowski }
 \end{center}

  \vspace{1cm}
 \noindent
 {\bf Abstract } In the present paper, we prove a sufficient
 condition of local regularity for suitable weak solutions to the
 Navier-Stokes equations having axial symmetry. Our condition is an
 axially symmetric analog of the so-called $L_{3,\infty}$-case in
 the
 general local regularity theory.

 \vspace {1cm}

\noindent {\bf 1991 Mathematical subject classification (Amer. Math.
Soc.)}: 35K, 76D.


\setcounter{equation}{0}
\section{Introduction  }
In the present paper, we address the problem of regularity for
axisymmetric solutions to the Navier-Stokes equations. In contrast
to many others, see, for example, \ci{L}, \ci{UY},\ci{LMNP},
\ci{NP},\ci{Po}, \ci{WZ}, \ci{Z1}, and \ci {Z2}, we study this
problem in local setting.

Our work is motivated by results of two different papers \ci{CL}
and \ci{ESS4}. To explain that, we need the following simple
notation. Let  $e_1$, $e_2$, $e_3$ be an orthogonal basis of the
Cartesian coordinates $x_1$, $x_2$, $x_3$ and $e_\varrho$,
$e_\varphi$, $e_3$ be an orthogonal basis  of the cylindrical
coordinates $\varrho$, $\varphi$, $x_3$ chosen so that
$$e_\varrho=\cos \varphi e_1+\sin \varphi e_2, \quad e_\varphi=-\sin\varphi
e_1+\cos\varphi e_2,\quad e_3=e_3.$$ Then, for any vector-valued
field $v$, we have  representations
$$v=v_ie_i=v_1e_1+v_2e_2+v_3e_3=v_\varrho e_\varrho+v_\varphi e_\varphi+v_3e_3.$$
The classical Navier-Stokes equations, which are invariant with
respect to transformation of coordinates, have the form \be\la{11}
\pa_t v +v\cdot \na v-\De v+\na p=0,\qquad \div v=0\ee and are
satisfied in some space-time domain. Here, as usual, $v$ and $p$
stand for the velocity field and the pressure field, respectively.

In our considerations, we always assume that $v_\varrho$,
$v_\varphi$, $v_3$, and $p$ are independent of the polar angle
$\varphi$. In \ci{CL}, Chae and Lee consider the Cauchy problem
for the Navier-Stokes equations under the above assumption on
axial symmetry. In addition to usual conditions on the initial
data, the authors of \ci{CL} assume that
  velocity field $v$ obeys
\be\la{12}\int\limits_0^Tdt\Big(\int\limits_{\mathbb R^3}
|v|^\ga d\varrho dx_3\Big)^\frac \al\ga<+\infty,\ee with
$1/\al+1/\ga\leq 1/2$, $2<\ga<+\infty$, $2< \al\leq +\infty$, and
prove the regularity of solutions to the Cauchy problem for
(\ref{11}) on  time interval $]0,T[$. In fact, they prove even
more, their statement is still true if $|v|$ is replaced with
$\sqrt{v_\varrho^2+v_\varphi^2}$. However, it remains unclear
whether or not the regularity takes place in the marginal case
$\ga=2$ and $\al=+\infty$. To our opinion, the case cannot be
treated by  methods developed in paper \ci{CL} because, in a
sense, it is an analog of the so-called $L_{3,\infty}$-case
studied in \ci{ESS4}. In turn,  $L_{3,\infty}$-case is marginal to
the so-called Ladyzhenskaya-Prodi-Serrin condition, see \ci{Pr},
\ci{Se}, \ci{L1a}, \ci{Str0}, \ci{Gi}, \ci{S7}, and \ci{S8}. It
seems quite reasonable to interpret the result of \ci{CL}, see
Theorem 3 there, as the Ladyzhenskaya-Prodi-Serrin condition for
axially symmetric problems.   To treat $L_{3,\infty}$-solutions in
generic setting, one needs  new technique based on backward
uniqueness for the heat operator with variable lower order terms.
In this, paper, we wish to extend this method to the axially
symmetric case.

To formulate our main result, we introduce the additional
notation:

$$\mathcal C(x_0,R)=\{x\in \mathbb R^3\,\,\|\,\, x=(x',x_3),
\,\,x'=(x_1,x_2),\,\,$$$$|x'-x'_0|<R,\,\,|x_3-x_{03}|<R\},\qquad\mathcal
C(R)=\mathcal C(0,R),\qquad \mathcal C=\mathcal C(1);
$$$$ z=(x,t),\qquad z_0=(x_0,t_0),\qquad Q(z_0,R)=\mathcal C(x_0,R)
\times ]t_0-R^2,t_0[, $$$$
  Q(R)=Q(0,R),\qquad Q=Q(1).$$

In local analysis, the most reasonable object to study is
so-called suitable weak solutions, introduced by
Caffarelli-Kohn-Nirenberg in their celebrated paper \ci{CKN}. We
are going to use a slightly simpler  definition of F.-H. Lin in
\ci{Li}
\begin{definition}\la{1d1} The pair $v$ and $p$ is called a
suitable weak solutions to the Navier-Stokes equations in $Q$ if
the following conditions are satisfied:
$$v\in L_{2,\infty}(Q)\cap W^{1,0}_2(Q),\qquad p \in L_\frac
32(Q);$$
$$v\,\mbox{and}\,p \,\mbox{satisfy the Navier-Stokes equtions in
the sense of distributions};$$ for a.a. $t\in ]-1,0[$, the local
energy inequality
$$\int\limits_\mathcal C\varphi(x,t)|v(x,t)|^2dx+2\int\limits^t_{-1}
\int\limits_\mathcal C\varphi |\na v|^2dxdt'\leq
\int\limits^t_{-1} \int\limits_\mathcal C\Big\{|v|^2(\De \varphi
+\pa_tv)+$$
$$+v\cdot\na \varphi(|v|^2+2p)\Big\}dxdt'$$
holds for all non-negative cut-off functions $\varphi\in
C^\infty_0(\mathbb R^3\times \mathbb R)$ vanishing in a
neighborhood of the parabolic boundary of $Q$.
\end{definition}
For discussions of the above definition, we refer the reader to
papers \ci{LS} and \ci{S8}.

Our main result is
\begin{theorem}\la{1t2} Let $v$ and $p$ be an axially symmetric suitable weak solution
to the Navier-Stokes equations in $Q$. Assume that
\be\la{13}\mathcal A_0= \ess \sup\limits_{-1\leq t\leq
0}\int\limits_\mathcal C\frac 1\varrho|v(x,t)|^2dx<+\infty.\ee
Then the point $(x,t)=(0,0)$ is a regular point of $v$, i.e.,
there exists $r\in ]0,1]$ such that $v$ is H\"older continuous  in
the closure of the cylinder $Q(r)$.
\end{theorem}

By $c$, we shall denote all generic constants that may vary from
one bound to others.

Our paper is organized as follows. In the second section, we
discuss known inequalities of the local regularity theory and
prove some useful facts about suitable weak solutions. The proof
of the main result is started in the third section with scaling
and blow up of our solution at a singular point. We also discuss
properties of the blow up velocity and the blow up pressure in
this section. In the fourth section, we prove some additional
differential properties of axially symmetric suitable weak
solutions. They are needed to establish a decay of the blow up
velocity at infinity. Finally, we end up with the proof of the
main theorem in the fifth section. Here, with the help of backward
uniqueness results for the heat operator with variable lower order
terms, we show that in fact our blow up velocity is trivial.

\noindent\textbf{Acknowledgement} The work is supported by the
Agreement on cooperation between Polish and Russian Academies of
Sciences signed in Warsaw, Dec. 27, 2002. The first author is
supported by the Alexander von Humboldt Foundation and by  the
RFFI grant 05-01-00941-a.

\setcounter{equation}{0}
\section{Preliminaries  }

In what follows, we are going to make use of the following scaling
invariant functionals:
$$A(z_0,r;v)=\ess \sup\limits_{t_0-r^2<t<t_0}\frac
1r\int\limits_{\mathcal C(x_0,r)}|v(x,t)|^2dx,\quad
C(z_0,r;v)=\frac 1{r^2}\int\limits_{Q(z_0,r)}| v|^3 dz, $$
$$E(z_0,r;v)=\frac
1r\int\limits_{Q(z_0,r)}|\na v|^2 dz,\qquad D(z_0,r;p)=\frac
1{r^2}\int\limits_{Q(z_0,r)}| p|^\frac 32 dz.$$

First, let us recall that, by the Navier-Stokes equations scaling,
$$v^\lb(x,t)=\lb v(\lb x,\lb^2t),\qquad p^\lb(x,t)=\lb^2 p(\lb
x,\lb^2t),$$ we may define suitable weak solutions to the
Navier-Stokes equations in $Q(z_0,R)$. So, if $v$ and $p$ form a
suitable weak solution to the Navier-Stokes equations in
$Q(z_0,R)$, then, for appropriate choice of the cut-off function
in the local energy inequality, we can reduce it to the following
invariant form $$A(z_0,R/2;v)+E(z_0,R/2;v)\leq c (C^\frac 23
(z_0,R;v)+$$\be\la{21}+C(z_0,R;v)+D(z_0,R;p)).\ee

We also need the so-called decay estimate for pressure
\be\la{22}D(z_0,r;p)\leq c \Big[\frac r{r_1}
D(z_0,{r_1};p)+\Big(\frac {r_1} r\Big)^2C(z_0,{r_1};v)\Big],\ee
which is valid for all $0<r\leq {r_1}\leq R$. The proof of the
latter estimate is given in \ci{S2}. Repeating arguments of Lemma
1.8 in \ci{S9}, we can prove

\begin{lemma}\la{2l1}Let $v$ and $p$ be a suitable weak solution
to the Navier-Stokes equations in $Q$ and let \be\la{23}
A_0=\sup\limits_{0<r<1}A(0,r;v)<+\infty.\ee Then, for any $r\in
]0,1/2[$, we have $$C^\frac 43(0,r;v)+D(0,r;p)+E(0,r;v)\leq c
\Big((A_0+1)r^\frac 12(D(0,1;p)+$$  \be\la{24}+
E(0,1;v))+A_0^4+A_0^2+A_0\Big).\ee
 \end{lemma}

Lemma \ref{2l1}, together with the invariance of our functionals
under the Navier-Stokes equations scaling and under the shift in
the direction of $x_3$, gives us:
\begin{lemma}\la{2l2}Under the conditions of Theorem \ref{1t2}, we have
\be\la{25}A(z_0,r;v)+C(z_0,r;v)+D(z_0,r;p)
+E(z_0,r;v)\leq\A<+\infty \ee for all $z_0=(x_0,0)$, $x_0=(0,b)$,
$|b|\leq 1/4$, and for all $0<r\leq 1/4$, where $\A$ depends on
$D(0,1;p)$, $E(0,1;v)$, and $\A_0$ only.
\end{lemma}

We say that the pair $v$ and $p$ is a suitable weak solution to
the Navier-Stokes equations in the space-time cylinder $\Om\times
]T_1,T_2[$,  if, for any $z_0=(x_0,t_0)$ with $x_0\in \Om$ and
$T_1<t\leq T_2$, the pair $v$ and $p$ is a suitable weak solution
to the Navier-Stokes equations in $Q(z_0,R)$ for some $R>0$.

Next, let us introduce the family of sets
$$\cl(R_1,R_2;a)=\{x\in \mathbb R^3\,\,\|\,\, R_1<|x'|<R_2,\, |x_3|<a\}.$$
Now, we would like to formulate and prove the following statement.
\begin{lemma}\la{2l3} Let $v$ and $p$ be a suitable
weak solution to the Navier-Stokes equations in the set
$\widehat{Q}=\cl(3/4,9/4;3/2)\times ]-(3/2)^2,0[$. Assume that
\be\la{26}\int\limits_{\widehat{Q}}|v(z)|^6dz\leq m<+\infty.\ee
Then, there exists a function $\Phi_0:\rb_+\times \rb_+\to\rb_+$,
nondecreasing in each variables, such  that
\be\la{27}|v(z)|+|\na v(z)|\leq \Phi_0(m,\A_*)<+\infty, \ee for
any $z\in \cl(1,2;1)\times ]-1,0[ $. Here,
$$\A_*=\int\limits_{\widehat{Q}}|p(z)|^\frac 32dz=
\int\limits_{-(3/2)^2}^0dt\int\limits_{\cl(3/4,9/4;3/2)}
|p|^\frac 32dx.$$
\end{lemma}
\textsc{Proof} First, we remark $Q(z_0,1/4)\subset\widehat{Q}$
for any
$z_0\in\cl(1,2;1)\times ]-1,0[$. It follows from (\ref{22}),
H\"older's inequality, and (\ref{26}) that
\be\la{28}D(z_0,r;p)\leq c \Big[\frac r{r_1}
D(z_0,{r_1};p)+\Big(\frac {r_1} r\Big)^2m^\frac 12r_1^\frac
12\Big],\ee which is valid for all $0<r\leq {r_1}\leq 1/4$. For
$\tau\in ]0,1[$, let us take $r=\tau^{k+1}/4$ and $r_1=\tau^{k}/4$
in (\ref{28}) and find
$$D(z_0,\tau^{k+1}/4;p)\leq c \tau\Big[
D(z_0,\tau^{k}/4;p)+ {m}^\frac 12{\tau}^{-3}\tau^{\frac k2}\Big]$$
for all non-negative integer numbers $k$. We can choose $\tau\in
]0,1[$ so small to provide
$$c\tau^\frac 34\leq 1$$
and conclude
$$D(z_0,\tau^{k+1}/4;p)\leq  \tau^\frac 14\Big[
D(z_0,\tau^{k}/4;p)+ {m}^\frac 12 {\tau}^{-3}\tau^{\frac
k2}\Big]$$ for all non-negative integer numbers $k$. The latter
inequality may be easily iterated. As a result, we have
$$D(z_0,\tau^{k+1}/4;p)\leq \tau^\frac {k+1}4\Big[D(z_0,1/4;p)
+ {m}^\frac 12  {\tau}^{-3}\sum\limits_{i=0}^k\tau^{\frac i
4}\Big]$$ for all non-negative integer numbers $k$. So,
$$C(z_0,\tau^{k+1}/4;v)+D(z_0,\tau^{k+1}/4;p)\leq
cm^\frac 12\tau^\frac{k+1}2+ \tau^\frac {k+1}4\Big[D(z_0,1/4;p)
$$$$+ {m}^\frac 12  {\tau}^{-3}(1-\tau^{\frac 14})^{-1}\Big]$$
$$\leq c\Big[m^\frac 12 \tau^\frac {k+1}2+ \tau^\frac
{k+1}4\Big(\A_*
+ {m}^\frac 12 {\tau}^{-3}(1-\tau^{\frac 14})^{-1}\Big)\Big]$$ for
all non-negative integer numbers $k$. Given $\varepsilon>0$, we
can find an integer number $k_0$ so that
$$c\Big[m^\frac 12 \tau^\frac {k_0+1}2+ \tau^\frac
{k_0+1}4\Big(\A_*
+ {m}^\frac 12 {\tau}^{-3}(1-\tau^{\frac 14})^{-1}\Big)\Big]\leq
\varepsilon.$$
But according to the so-called $\varepsilon$-regularity theory,
see, for example, \ci{LS}, \ci{ESS4}, and \ci{S8}, the latter
implies two bounds:
$$| v(z_0)|\leq \frac c{r_0}\quad\mbox{and}\quad |\na v(z_0)|\leq \frac c{r_0^2},$$
where $r_0=\tau^{(k_0+1)}/4$. Lemma \ref{2l3} is proved.

The last preliminary statement is as follows.

\begin{lemma}  \la{2l4} Assume that all conditions of
Theorem \ref{1t2} hold. Then \be\la{29}\int\limits_\mathcal C\frac
1\varrho|v(x,t)|^2dx\leq\mathcal A_0\ee for all $ t\in
]-1,0[$.\end{lemma} \textsc{Proof} It easy to derive the following
estimate
$$\int\limits_Q\pa_t v\cdot wdz\leq \A_1\Big(
\int\limits_Q|\na w|^3dz\Big)^\frac 13$$ for any $C^\infty_0(Q)$.
Here, a constant $\A_1$ depends on $C(0,1;v)$, $E(0,1;v)$, and
$D(0,1;p)$ only. So, $v$ has the first derivative in to $t$ in the
space $$L_\frac 32(-1,0; W^{-1}_\frac 32(\C)).$$ In turn, the
latter, together with boundedness of the energy, implies weak
continuity in time in the following sense: the function
$$t\to \int\limits_\C v(x,t)\cdot w(x)dx$$ is continuous on
$[-1,0]$ for any $w\in L_2(\C)$. Now, the statement of the lemma
follows from the weak lower  semicontinuity  of the functional
$$w\in L_2(\C)\to \int\limits_\C\frac 1\varrho|w(x)|^2dx.$$
Lemma \ref{2l4} is proved.

\setcounter{equation}{0}
\section{Scaling and Blow Up }

Here, we are starting with the proof of Theorem \ref{1t2}. Assume
that the statement of this theorem is false. Then, according to  the
local regularity theory for the Navier-Stokes equations, there exist
an absolute positive constant $\varepsilon$ and a sequence
$\{R_k\}^\infty_{k=1}$ such that $R_k\to 0$ as $k\to +\infty$ and
\be\la{31}\frac 1{R_k^2}\int\limits_{Q(R_k)}|v|^3dz\geq
\varepsilon>0\ee for all $k\in \mathbb N$.

Next, we scale $v$ and $p$ in the following way
$$u^k(y,s)=R_kv(R_ky,R_k^2s),\quad
q^k(y,s)=R^2_kp(R_ky,R_k^2s),$$ where $e=(y,s)\in Q(1/R_k)$.
Functions $u^k$ and $q^k$ are extended by zero to the whole
space-time $\rb^3\times \mathbb R$.

Now let us fix numbers $a$ and $b$ in $\rb$ so that $a>0$.
Let
$$x_k^b=(0,bR_k),\quad y^b=(0,b),\quad z_k^b=(x_k^b,0),\quad
e^b=(y^b,0).$$ Obviously, for sufficiently large $k$,
$$|b|R_k<1/4,\qquad aR_k<1/4,$$
by Lemma \ref{2l2}, the following estimates are valid:
$$C( z_k^b,aR_k;v)=C(e^b,a;u^k)\leq \A,$$
$$E( z_k^b,aR_k;v)=E(e^b,a;u^k)\leq \A,$$
\be\la{32}A( z_k^b,aR_k;v)=A(e^b,a;u^k)\leq \A,\ee
$$D( z_k^b,aR_k;p)=D(e^b,a;q^k)\leq \A$$
for all $k\geq k_0(a,b)$.

First, let $b$ be equal to zero. In this particular case, we can
produce three estimates. The first of them is well known in the
Navier-Stokes theory and it is but a consequence of multiplicative
inequalities \be\la{33}\frac 1{a^\frac
52}\int\limits_{Q(a)}|u^k|^\frac {10}3de\leq c(\A).\ee The second
estimate follows from the Navier-Stokes equations, written for $u^k$
and $q^k$ in the weak form, and from (\ref{32}):
$$\int\limits_{Q(a)}\pa_tu^k\cdot wde\leq c(a,\A)\Big(\int\limits_{Q(a)}
|\na w|^3de\Big)^\frac 13$$ for all $w\in C^\infty_0(Q(a))$. Hence,
\be\la{34}\pa_tu^k\,\mbox{is bounded in} \, L_\frac
32(-a^2,0;W_\frac 32^{-1}(\C(a))).\ee The third estimate is coming
from our main condition (\ref{13}) and has the form \be\la{35}\ess
\sup\limits_{-(aR_k)^2\leq t\leq 0}\int\limits_{\mathcal
C(aR_k)}\frac {|v(x,t)|^2}{|x'|}dx=\ess \sup\limits_{-a^2\leq s\leq
0}\int\limits_{\mathcal C(a)}\frac
{|u^k(y,t)|^2}{|y'|}dy\leq\A_0.\ee

Now, making use of the diagonal process for extending space-time
cylinders $Q(a)$ and known compactness arguments, we can select
subsequences (still denoted by $u^k$ and $q^k$) such that, for each
$a>0$,
$$u^k{\rightharpoondown }\, u\qquad \mbox{in}\,\,W^{1,0}_2(Q(a)),$$
$$u^k\stackrel{\star}{\rightharpoondown} u\qquad \mbox{in}\,\,L_{2,\infty}(Q(a)),$$
\be\la{36}u^k\rightarrow u\qquad \mbox{in}\,\,L_3(Q(a)),\ee
$$q^k{\rightharpoondown } \,q\qquad \mbox{in}\,\,L_\frac 32(Q(a)).$$
The aim of our further considerations is to describe properties of
limit functions $u$ and $q$ called the blow up velocity and blow up
pressure, respectively. They are defined on $\rb^3\times \rb_-$,
where $\rb_-=\{s\in \rb\,\,\|\,\,s\leq 0\}$. For each $a>0$, the
pair $u$ and $q$ is a suitable weak solution to the Navier-Stokes
equations in $Q(a)$. From (\ref{32}) and (\ref{36}), it follows that
the limit functions obey the inequalities:
$$C(e^b,a;u)\leq \A,$$
$$A(e^b,a;u)\leq \A,$$
\be\la{37}E(e^b,a;u)\leq \A,\ee
$$D(e^b,a;q)\leq \A$$
for all $b\in\rb$ and for all $0<a\in \rb$. Moreover, we can
derive from (\ref{36}), (\ref{35}), and (\ref{31}) two additional
estimates: \be\la{38}\ess \sup\limits_{-\infty< s\leq
0}\int\limits_{\rb^3}\frac {|u(y,t)|^2}{|y'|}dy\leq\A_0\ee and
\be\la{39}\frac
1{R_k^2}\int\limits_{Q(R_k)}|v|^3dz=\int\limits_{Q}|u^k|^3de\rightarrow
\int\limits_{Q}|u|^3de\geq\varepsilon.\ee

According to (\ref{39}), the blow up velocity $u$ is a non-trivial
solution to the Navier-Stokes equations in $\rb^3\times\rb_-$. But
we are going to show that in fact $u\equiv 0$. This would contradict
with (\ref{39}) and prove Theorem \ref{1t2}.

Obviously, the blow up velocity field $u$ is axially symmetric and,
by Caffarerrli-Kohn-Nirenberg type results, all point $y'\neq 0$ are
regular which make it possible to conclude that all spatial
derivatives of $u$ are H\"older continuous in a vicinity of each
point with $y'\neq 0$.

We  can also make  use of the local regularity theory for Stokes
system, see \ci{S7} and \ci{S8}. According to it and by known
multiplicative inequality, we have
$$\|\pa_t u^k\|_{L_{\frac 98, \frac 32}(Q(a/2))}+
\|\na^2 u^k\|_{L_{\frac 98, \frac 32}(Q(a/2))}+\|\na
q^k\|_{L_{\frac 98, \frac 32}(Q(a/2))}\leq$$
$$\leq c(a)\Big[\|u^k\cdot\na u^k\|_{L_{\frac 98, \frac 32}(Q(a))}
+\|u^k\|_{W^{1,0}_2(Q(a))}+\|q^k\|_{L_\frac 32(Q(a))}\Big]\leq$$
$$\leq c(a)\Big[\|u^k\|^\frac 23_{L_{2,\infty}(Q(a))}
\|u^k\|_{W^{1,0}_2(Q(a))}^\frac 13+...\Big]\leq c(a,\A).$$ The
latter estimate shows that we can select a subsequence (still
denoted by $u^k$) such that, for any $a>1$, \be\la {310}
u^k\rightarrow\,u\qquad\mbox{in}\,\, C([-1,0];L_\frac 98(\C
(a))).\ee (\ref{310}) can be exploited in the following way. For any
fixed positive numbers $r_1$, $r_2$, and $h$, we have
$$\Big(\int\limits_{\cl(r_1,r_2;h)}|u(y,0)|^\frac 98dy\Big)^\frac 89\leq
\Big(\int\limits_{\cl(r_1,r_2;h)}|u^k(y,0)-u(y,0)|^\frac
98dy\Big)^\frac 89+$$$$
+\Big(\int\limits_{\cl(r_1,r_2;h)}|u^k(y,0)|^\frac 98dy\Big)^\frac
89=\al_k+\bt_k.$$ By (\ref{310}),
$$\al_k\rightarrow 0$$
as $k\to+\infty$. To evaluate $\bt_k$, we make use of the inverse
scaling and H\"older's inequality
$$\bt_k=\Big({R^{-\frac {15}8}_k}\int\limits_
{\cl(R_kr_1,R_kr_2;R_kh)}|v(x,0)|^\frac 98dx\Big)^\frac 89\leq $$
$$\leq c(r_1,r_2,h)\Big(\frac
1{R_k}\int\limits_{\cl(R_kr_1,R_kr_2;R_kh)}|v(x,0)|^2dx\Big)^\frac
12\leq $$
$$\leq c(r_1,r_2,h)\Big(\int\limits_{\cl(R_kr_1,R_kr_2;R_kh)}
\frac {|v(x,0)|^2}{|x'|}dx\Big)^\frac 12.$$ Now, it remains to apply
Lemma \ref{2l4} at $t=0$ and absolute continuity of Lebesgue's
integral and conclude that
$$\bt_k\rightarrow\,0$$
as $k\to+\infty$. This implies the identity
$$\int\limits_{\cl(r_1,r_2;h)}|u(y,0)|^\frac 98dy=0$$ for all
positive numbers $r_1$, $r_2$, and $h$. So, we can state that
\be\la{311} u(\cdot,0)=0\qquad\mbox{in}\quad\rb^3.\ee

\setcounter{equation}{0}
\section{Estimates of Axially Symmetric Solutions  }

The main result of this section is going to be as follows.
\begin{pro}\la{4p1} Let $V$ and $P$ be a sufficiently
smooth axially symmetric solution to the Navier-Stokes equations
in $\widetilde{Q}=\widetilde{\cl}\times ]-2^2,0[$, where
$\wcl=\cl(1/4,3;2)$. Then, there exists a non-decreasing function
$\Phi:\rb_+\to\rb_+$ such that
\be\la{41}\sup\limits_{z\in\cl(1,2;1)\times
]-1,0[}\Big(|V(z)|+|\na V(z)|\Big)\leq \Phi(\A_2),\ee where
$$\A_2=\sup\limits_{-2^2<t<0}\int\limits_{\wcl}|V(x,t)|^2dx+\int\limits_
{\wq}\Big(|\na V|^2+|V|^3+|P|^\frac 32\Big)dz.$$
\end{pro}

To prove the above proposition, we need
\begin{lemma}\la{4l2} Under assumptions of Proposition \ref{4p1},
there exists a function $\Phi_1:\rb_+\times\rb_+\to\rb_+$,
non-decreasing in each variable, such that \be\la{42}\sup
\limits_{-(7/4)^2<t<0}\int\limits_{\wcl_1}|V^a(x,t)|^q dx\leq
\Phi_1(q,\A_2),\qquad 1\leq q +\infty.\ee Here,
$V^a=(V_\varrho,V_3)$, $|V^a|=\sqrt{|V_\varrho|^2+|V_3|^2}$,
$\wcl_1=\cl(5/16,11/4;7/4)$, and $\widetilde{Q}_1=\wcl_1\times
]-(7/4)^2,0[$.\end{lemma} \textsc{Proof} Let us denote by $\om$
the vorticity of $v$, i.e., $\om=\na \wedge v$. For
$\chi=\om_\varphi$, $V_\varrho$, and $V_3$, we have the following
identities: \be\la{43} V_{\varrho,\varrho}+V_{3,3}=-\frac 1\varrho
V_\varrho,\ee \be\la{44}V_{\varrho,3}-V_{3,\varrho}=\chi,\ee
$$\pa_t\chi+V_\varrho\chi_{,\varrho}+V_3\chi_{,3}-\frac
1\varrho\chi V_\varrho-\Big(\chi_{,\varrho\varrho}+\chi_{,33}
+\frac 1\varrho\chi_{,\varrho}-\frac
1{\varrho^2}\chi\Big)=$$\be\la{45}=\frac 2\varrho V_\varphi
V_{\varphi,3},\ee where we have used the notion
$$f_{,\varrho}=\frac {\pa f}{\pa \varrho},\qquad
f_{,3}=\frac {\pa f}{\pa x_3}.$$ Next, we let
$\widetilde{\chi}=\chi\psi$,
 $\widetilde{V}=V^a\psi$, $\widetilde{V}_\varrho=V_\varrho\psi$, and
$\widetilde{V}_3=V_3\psi$, where a non-negative smooth and axially
symmetric cut-off function $\psi$ vanishes in a neighborhood of
the parabolic boundary of $\wq$ and is equal to 1 in $\wq_1$. For
$\widetilde{\chi}$, $\widetilde{V}_\varrho$, and
$\widetilde{V}_3$, we have
\be\la{46}\widetilde{V}_{\varrho,\varrho}+\widetilde{V}_{3,3}=-\frac
1\varrho
\widetilde{V}_\varrho+V_\varrho\psi_{,\varrho}+V_3\psi_{,3},\ee
\be\la{47}\widetilde{V}_{\varrho,3}-\widetilde{V}_{3,\varrho}=
\widetilde{\chi}+V_\varrho\psi_{,3}-V_3\psi_{,\varrho},\ee
$$\pa_t\widetilde{\chi}+V_\varrho\widetilde{\chi}_{,\varrho}+
V_3\widetilde{\chi}_{,3}-\frac 1\varrho
V_\varrho\widetilde{\chi}-\Big(\widetilde{\chi}_{,\varrho\varrho}+\widetilde{\chi}_{,33}
+\frac 1\varrho\widetilde{\chi}_{,\varrho}-\frac
1{\varrho^2}\widetilde{\chi}\Big)=$$\be\la{48}=J_1+J_2+J_3, \ee
where $$J_1=\frac 2\varrho V_\varphi
V_{\varphi,3}\psi,$$$$J_2=\chi\Big(\pa_t\psi
-\psi_{,\varrho\varrho}-\psi_{,33}-\frac
1\varrho\psi_{,\varrho}\Big)-2\Big(\chi_{,\varrho}
\psi_{,\varrho}+\chi_{,3}\psi_{,3}\Big),$$
$$J_3=\chi\Big(V_\varrho\psi_{,\varrho}+V_3\psi_{,3}\Big).$$

Now, we multiply (\ref{48}) by $\widetilde{\chi}\varrho^{-2}$ and
integrate the product by parts over $\wcl$ $$\frac 12\pa_t
\int\limits_{\wcl}\Big|\frac
{\widetilde{\chi}}\varrho\Big|^2dx+\int\limits_{\wcl}\Big(\Big|\Big
(\frac {\widetilde{\chi}}\varrho\Big)_{,\varrho}\Big|^2+\Big|\Big
(\frac {\widetilde{\chi}}\varrho\Big)_{,3}\Big|^2\Big)dx=$$
\be\la{49}=\int\limits_{\wcl}J_1\frac
{\widetilde{\chi}}{\varrho^2}dx+\int\limits_{\wcl}J_2\frac
{\widetilde{\chi}}{\varrho^2}dx+\int\limits_{\wcl}J_3\frac
{\widetilde{\chi}}{\varrho^2}dx.\ee Our aim is to evaluate the
right hand side of (\ref{49}). We start with the first term there:
$$\int\limits_{\wcl}J_1\frac
{\widetilde{\chi}}{\varrho^2}dx=-\int\limits_{\wcl}\frac
{V_\varphi^2}{\varrho^2}\Big(\frac
{\widetilde{\chi}}{\varrho}\Big)_{,3}\psi
dx-\int\limits_{\wcl}\frac {V_\varphi^2}{\varrho^2}\frac
{\widetilde{\chi}}{\varrho}\psi_{,3} dx\leq$$
$$\leq c\Big(\int\limits_{\wcl}\frac
{|V_\varphi|^4}{\varrho^4}dx\Big)^\frac 12
\Big(\int\limits_{\wcl}\Big|\Big (\frac
{\widetilde{\chi}}\varrho\Big)_{,3}\Big|^2dx+\int\limits_{\wcl}\Big|\frac
{\widetilde{\chi}}\varrho\Big|^2dx\Big)^\frac 12,$$
where the notion $\na_af=(f_{,\varrho},f_{,3})$ has been used. To
estimate the first multiplier of the right hand of the latter
inequality, we are going to exploit two-dimensional feature of our
axially symmetric problem in the following way. So, by
Ladyzhenskaya's inequality,
$$\int\limits_{\wcl}
{|V_\varphi|^4}dx\leq c\int\limits_{-2}^2\int\limits_{1/4}^3
|V_\varphi|^4d\varrho dx_3\leq$$$$\leq
c\int\limits_{-2}^2\int\limits_{1/4}^3|V_\varphi|^2d\varrho dx_3
\int\limits_{-2}^2\int\limits_{1/4}^3\Big(|V_\varphi|^2
+|\na_aV_\varphi|^2\Big)d\varrho dx_3\leq$$
$$\leq c\int\limits_{\wcl}|V|^2dx\int\limits_{\wcl}
\Big(|V|^2 +|\na V|^2\Big)dx\leq c\A_2\int\limits_{\wcl}
\Big(|V|^2 +|\na V|^2\Big)dx.$$ Thus, we find the first estimate:
$$\int\limits_{\wcl}J_1\frac
{\widetilde{\chi}}{\varrho^2}dx \leq c\A_2^\frac 12 \Big(
\int\limits_{\wcl} \Big(|V|^2 +|\na V|^2\Big)dx\Big)^\frac
12\times$$\be\la{410}\times\Big(\int\limits_{\wcl}\Big|\na_a\Big
(\frac
{\widetilde{\chi}}\varrho\Big)\Big|^2dx+\int\limits_{\wcl}\Big|\frac
{\widetilde{\chi}}\varrho\Big|^2dx\Big)^\frac 12.\ee

For the second term, we have \be\la{411}\int\limits_{\wcl}J_2\frac
{\widetilde{\chi}}{\varrho^2}dx\leq c
\int\limits_{\wcl}|\chi|^2dx.\ee The third term is estimated in
slightly different way
$$\int\limits_{\wcl}J_3\frac
{\widetilde{\chi}}{\varrho^2}dx=\int\limits_{\wcl}\frac
{\chi\widetilde{\chi}}{\varrho^2}
\Big(V_\varrho\psi_{,\varrho}+V_3\psi_{,3}\Big)dx\leq
$$$$\leq c\Big(\int\limits_{\wcl}|\chi|^2dx\Big)^\frac 12
\Big(\int\limits_{\wcl}|V^a\cdot \na_a\psi|^2\Big|\frac
{\widetilde{\chi}}\varrho\Big|^2dx\Big)^\frac 12\leq$$
$$\leq c\int\limits_{\wcl}|\chi|^2dx+c\Big(\int\limits_{\wcl}|
V^a\cdot \na_a\psi|^4dx\Big)^\frac
12\Big(\int\limits_{\wcl}\Big|\frac
{\widetilde{\chi}}\varrho\Big|^4dx\Big)^\frac 12,$$ where we let
$V^a\cdot \na_a\psi=V_\varrho\psi_{,\varrho}+V_3\psi_{,3}$. To
estimate the last term on the right hand side of the latter
relation, we exploit Ladyzhenskaya's inequality once more. So, we
have
$$\int\limits_{\wcl}\Big|\frac
{\widetilde{\chi}}\varrho\Big|^4dx\leq
c\int\limits_{-2}^2\int\limits_{1/4}^3\Big|\frac
{\widetilde{\chi}}\varrho\Big|^4d\varrho dx_3\leq$$
$$\leq c\int\limits_{-2}^2\int\limits_{1/4}^3\Big|
\na_a\Big(\frac {\widetilde{\chi}}\varrho\Big)\Big|^2d\varrho
dx_3\int\limits_{-2}^2\int\limits_{1/4}^3\Big|\frac
{\widetilde{\chi}}\varrho\Big|^2d\varrho dx_3\leq$$
$$\leq c\int\limits_{\wcl}\Big|
\na_a\Big(\frac {\widetilde{\chi}}\varrho\Big)\Big|^2dx
\int\limits_{\wcl}\Big|\frac {\widetilde{\chi}}\varrho\Big|^2dx$$
and, in the same way,
$$\int\limits_{\wcl}|
V^a\cdot \na_a\psi|^4dx\leq c\int\limits_{\wcl}\Big|
\na_a\Big(V^a\cdot \na_a\psi\Big)\Big|^2dx \int\limits_{\wcl}|
V^a\cdot \na_a\psi|^2dx.$$ As a result, we find
$$\int\limits_{\wcl}J_3\frac
{\widetilde{\chi}}{\varrho^2}dx\leq c\int\limits_{\wcl}|\chi|^2dx+
c\Big(\int\limits_{\wcl}|V_a|^2dx+
\int\limits_{\wcl}|\na_aV^a|^2dx\Big)^\frac 12\times$$$$\times
\Big(\int\limits_{\wcl}|V_a|^2dx\Big)^\frac
12\Big(\int\limits_{\wcl}\Big| \na_a\Big(\frac
{\widetilde{\chi}}\varrho\Big)\Big|^2dx\Big)^\frac 12
\Big(\int\limits_{\wcl}\Big|\frac
{\widetilde{\chi}}\varrho\Big|^2dx\Big)^\frac 12\leq$$
\be\la{412}\leq c\int\limits_{\wcl}|\na V|^2dx+c\Big(\A_2
+\A_2^\frac 12\Big(\int\limits_{\wcl}|\na V|^2dx\Big)^\frac
12\Big)\times\ee
$$\times\Big(\int\limits_{\wcl}\Big| \na_a\Big(\frac
{\widetilde{\chi}}\varrho\Big)\Big|^2dx\Big)^\frac 12
\Big(\int\limits_{\wcl}\Big|\frac
{\widetilde{\chi}}\varrho\Big|^2dx\Big)^\frac 12.$$ Combining
estimates (\ref{49})-(\ref{412}) and applying Young's inequality,
we arrive at the final inequality
$$\pa_t
\int\limits_{\wcl}\Big|\frac {\widetilde{\chi}}\varrho\Big|^2dx+
\int\limits_{\wcl}\Big| \na_a\Big(\frac
{\widetilde{\chi}}\varrho\Big)\Big|^2dx\leq c\int\limits_{\wcl}|\na
V|^2dx+$$ \be\la{413}+\Big(\A_2^2 +\A_2\int\limits_{\wcl}|\na
V|^2dx\Big)\Big(\int\limits_{\wcl}\Big|\frac
{\widetilde{\chi}}\varrho\Big|^2dx+1\Big).\ee Estimate (\ref{413})
implies
$$\|\widetilde{\chi}\|_{L_{2,\infty}(\widetilde{Q})}\leq \Phi_3(\A_2).$$
According to (\ref{46}) and (\ref{47}), one may conclude
$$\int\limits_{-2}^2\int\limits_{1/4}^3\Big |\na_a
\widetilde{V}\Big|^2d\varrho dx_3\leq
c\int\limits_{-2}^2\int\limits_{1/4}^3\Big(|\widetilde{\chi}|^2+|V^a|^2
\Big)d\varrho dx_3\leq\Phi_3(\A_2)$$ and thus
$$\int\limits_{-2}^2\int\limits_{1/4}^3\Big |
\widetilde{V}(x,t)\Big|^qd\varrho dx_3\leq\Phi_4(q,\A_2)$$ for all
$t\in ]-2^2,0[$. Now, (\ref{42}) immediately follows  from the
latter inequality. Lemma \ref{4l2} is proved.

The second counterpart of the proof of Proposition \ref{4p1} is the
following statement.
\begin{lemma}\la{4l3} Under assumptions of Proposition \ref{4p1},
there exists a non decreasing function $\Phi_5:\rb_+\to\rb_+$ such
that \be\la{414}\int\limits_{\wq_2}|V_\varphi|^6dz\leq
\Phi_5(\A_2),\ee where $\wq_2=\wcl_2\times ]-(3/2)^2,0[$ and
$\wcl_2=\cl(3/8,5/2;3/2)$.\end{lemma} \textsc{Proof} We know that
$V_\varphi$ satisfies the equation
$$\pa_tV_\varphi+V_\varrho
V_{\varphi,\varrho}+V_3V_{\varphi,3}+\frac 1\varrho V_\varrho
V_\varphi$$\be\la{415}-\Big(V_{\varphi,\varrho\varrho}+
V_{\varphi,33}+\frac 1\varrho V_{\varphi,\varrho} -\frac
1{\varrho^2}V_\varphi\Big)=0.\ee

We fix a non-negative smooth and axially symmetric cut-off function
$\psi$ vanishing in a neighborhood of the parabolic boundary of $
\wq_1$ and being equal to 1 in $\wq_2$. Then, for
$\widetilde{\al}=V_\varphi\psi\varrho$,  we have the following
identity
$$\pa_t\wal+V_\varrho\wal_{,\varrho}+V_3\wal_{,3}
-\Big(\wal_{,\varrho\varrho}+\wal_{,33}+\frac
1\varrho\wal_{,\varrho}\Big)+\frac 2\varrho\wal_{,\varrho}=$$
$$=\al \Big(\pa_t\psi+V_\varrho\psi_{,\varrho}+
V_3\psi_{,3}\Big)-\Big(2\al_{,\varrho}\psi_{,\varrho}+2\al_{,3}
\psi_{,3}+\al\psi_{,\varrho\varrho}+\al\psi_{,33}\Big)+\frac
1\varrho\al\psi_{,\varrho},$$
where $\al=V_\varphi\varrho$.

Then, we multiply the latter identity by $\wal|\wal|^2$ and
integrate the product by parts over $\wcl_1$ \be\la{416} \frac
14\pa_t\int\limits_{\wcl_1}|\wal|^4dx+\frac
34\int\limits_{\wcl_1}|\na_a(|\wal|^2)|^2dx=J_1+J_2,\ee where
$$J_1=\int\limits_{\wcl_1}\al\wal|\wal|^2\Big(
V_\varrho\psi_{,\varrho}+V_3\psi_{,3}\Big)dx$$ and
$$J_2=\int\limits_{\wcl_1}\wal|\wal|^2\Big(\al\pa_t\psi-
2\al_{,\varrho}\psi_{,\varrho}-2\al_{,3}\psi_{,3}
-\al\psi_{,\varrho\varrho}-\al\psi_{,33}+\frac
1\varrho\al\psi_{,\varrho}\Big)dx.$$ We let $\bt=|\wal|^2$, then
$|\bt|^\frac {10}3=|\wal|^\frac {20}3$ and \be\la{417}
\int\limits_{\wcl_1}|\bt|^\frac {10}3dx\leq
c\Big(\int\limits_{\wcl_1}|\bt|^2dx\Big)^\frac
23\int\limits_{\wcl_1}|\na\bt|^2dx.\ee

We start with $J_1$, setting $\A_3=\|V^a\|_{L_{4,\infty}(\wq_1)}$.
By H\"older's inequality and by  multiplicative inequality
(\ref{417}),
$$J_1\leq c\Big(\int\limits_{\wcl_1}|\wal|^\frac {20}3
dx\Big)^\frac 9{20}\Big(\int\limits_{\wcl_1}|\al|^\frac {20}{11}
|V^a|^\frac {20}{11}dx\Big)^\frac {11}{20}\leq$$
$$\leq c\Big(\int\limits_{\wcl_1}|\bt|^\frac {10}3
dx\Big)^\frac 9{20}\Big(\int\limits_{\wcl_1}|\al|^\frac {10}3
dx\Big)^\frac 3{10}\Big(\int\limits_{\wcl_1}|V^a|^4dx\Big)^\frac
14\leq$$$$\leq c\Big(\int\limits_{\wcl_1}|\bt|^\frac {10}3
dx\Big)^\frac 9{20}\Big(\int\limits_{\wcl_1}|V|^2dx\Big)^\frac
15\Big(\int\limits_{\wcl_1}|V|^2dx+\int\limits_{\wcl_1}|\na
V|^2dx\Big)^\frac 3{10}\A_3\leq$$$$\leq
c\Big(\int\limits_{\wcl_1}|\bt|^\frac {10}3 dx\Big)^\frac
9{20}\A_2^\frac 15\Big(\int\limits_{\wcl_1}|\na
V|^2dx+\A_2\Big)^\frac 3{10}\A_3\leq$$$$\leq
c\Big(\int\limits_{\wcl_1}|\bt|^2dx\Big)^\frac
3{10}\Big(\int\limits_{\wcl_1}|\na\bt|^2dx\Big)^\frac
9{20}\A_2^\frac 15\Big(\int\limits_{\wcl_1}|\na
V|^2dx+\A_2\Big)^\frac 3{10}\A_3.$$ Term $J_2$ is estimated in the
same way:
$$J_2\leq c\Big(\int\limits_{\wcl_1}|\wal|^\frac {20}3
dx\Big)^\frac 9{20}\Big(\int\limits_{\wcl_1}(|\al|+|\al_{,\varrho}|
+|\al_{,3}|)^\frac {20}{11} dx\Big)^\frac {11}{20}\leq$$
$$\leq
c\Big(\int\limits_{\wcl_1}|\bt|^2dx\Big)^\frac
3{10}\Big(\int\limits_{\wcl_1}|\na\bt|^2dx\Big)^\frac
9{20}\Big(\int\limits_{\wcl_1}|\na V|^2dx+\A_2\Big)^\frac 12.$$
Now, making use of Young's inequality, we derive from (\ref{417})
and from  two latter estimates the main inequality
$$\pa_t\int\limits_{\wcl_1}|\wal|^4dx+
\int\limits_{\wcl_1}|\na_a(|\wal|^2)|^2dx\leq $$$$\leq
c\Big(\int\limits_{\wcl_1}|\bt|^2dx\Big)^\frac 6{11} \A_2^\frac
4{11}\Big(\int\limits_{\wcl_1}|\na V|^2dx+\A_2\Big)^\frac
6{11}\A_3^\frac
{20}{11}$$$$+c\Big(\int\limits_{\wcl_1}|\bt|^2dx\Big)^\frac 6{11}
\Big(\int\limits_{\wcl_1}|\na V|^2dx+\A_2\Big)^\frac {10}{11}\leq$$
$$\leq c\int\limits_{\wcl_1}|\bt|^2dx
\Big(\int\limits_{\wcl_1}|\na V|^2dx+\A_2\Big)
+c(\A_2\A_3^5)^\frac 4{11}+c\Big(\int\limits_{\wcl_1}|\na
V|^2dx+\A_2\Big)^\frac 4{11}.$$ It, together with the statement of
Lemma \ref{4l2} at $q=4$, implies \be\la{418}\sup_{-(7/4)^2\leq
t\leq 0}
\int\limits_{\wcl_1}|\bt(x,t)|^2dx+\int\limits_{\wq_1}|\na\bt|^2dz
\leq \Phi_5(\A_2).\ee So, (\ref{414}) follows from (\ref{417}) and
(\ref{418}). Lemma \ref{4l3} is proved.

From Lemmata \ref{4l2} and \ref{4l3}, we find
\begin{cor}\la{4c4} Under assumptions of Proposition \ref{4p1},
there exists a non-decreasing function $\Phi_6:\rb_+\to\rb_+$ such
that \be\la{419}\int\limits_{\wq_2}|V|^6dz\leq
\Phi_6(\A_2).\ee\end{cor} \textsc{Proof of Proposition \ref{4p1}}
Applying Corollary \ref{4c4} and Lemma \ref{2l3}, we end up with the
proof of Proposition \ref{4p1}. Proposition \ref{4p1} is proved.

\setcounter{equation}{0}
\section{Proof of Theorem \ref{1t2}  }

Given $R>1$, let us consider the following space-time cylinder
$$\wq^b_R=\wcl^b_R\times ]-(2R)^2,0[,$$
where $b\in\rb$ and
$$\wcl^b_R=\wcl_R+be_3,\qquad \wcl_R=\cl(R/4,3R;2R).$$
Now, we scale our blow up functions $u$ and $q$ in the following way
$$u^R(x,t)=Ru(Rx+be_3,R^2t),\qquad q^R(x,t)=R^2q(Rx+be_3,R^2t)$$
for $z=(x,t)\in \wq$.

Functions $u^R$ and $q^R$ are axially symmetric and, as it was
explained before, sufficiently smooth to apply Proposition
\ref{4p1}. According to that, we have
$$\sup_{z\in\wq_0}\Big\{|u^R(z)|+|\na u^R(z)|\Big\}\leq \Phi(\A_2),$$
where $\wq_0=\cl(1,2;1)$ and
$$\A_2=\sup_{-2^2\leq t\leq 0}\int\limits_{\wcl}|u^R(x,t)|^2dx
+\int\limits_{\wq}\Big(|\na u^R|^2+|u^R|^3+|q^R|^\frac 32\Big)dz.$$
Then, we make the inverse change of variables. As a result, we find
$$\sup_{(y,s)\in Q^b_R}\Big\{R|u(y,s)|+R^2|\na u(y,s)|\Big\}\leq \Phi
(\widetilde{\A}_{2R}),$$ where $Q^b_R=\wcl^b_{0R}\times ]-R^2,0[$,
$\wcl^b_{0R}=be_3+\cl_{0R}$, $\cl_{0R}=\cl(R,2R;R)$, and
$$\widetilde{\A}_{2R}=\sup\limits_{-(2R)^2\leq s\leq 0}\frac 1R
\int\limits_{-2R+b}^{2R+b}dy_3\int\limits_{R/4<|y'|<3R}|u(y,s)|^2dy'+$$
$$+\frac 1R
\int\limits_{-(2R)^2}^0ds\int\limits_{-2R+b}^{2R+b}dy_3
\int\limits_{R/4<|y'|<3R}|\na u(y,s)|^2dy'+$$
$$+\frac 1{R^2}
\int\limits_{-(2R)^2}^0ds\int\limits_{-2R+b}^{2R+b}dy_3
\int\limits_{R/4<|y'|<3R}\Big(| u(y,s)|^3+|q(y,s)|^\frac
32\Big)dy'\leq$$
$$\leq
c\Big(A(e^b,3R;u)+E(e^b,3R;u)+C(e^b,3R;u)+D(e^b,3R;q)\Big)\leq
c\A,$$ $e^b=(y^b,0)$ and $y^b=(0,b)$. So, assuming that $|y'|>20$,
we can derive from the latter estimates
\be\la{51}|y'||u(y',b,s)|+|y'|^2|\na u(y',b,s)|\leq \Phi(c\A)\ee
for any $b\in \rb$, for any $|y'|>20$, and for any $s\in [-20,0]$.
It follows directly from (\ref{51}) that: \be\la{52}|u(y,s)|+|\na
u(y,s)|\leq c \Phi(c\A)=c(\A)\ee for any $|y'|>20$ and for any
$s\in [-20,0]$.

Now, we consider the vorticity $\om(u)=\na \wedge u$. It satisfies
the vorticity equation
$$\pa_t\om-\De\om=\om\cdot\na u-u\cdot \na\om,$$
which, together with \ref{52}, implies
\be\la{53}|\pa_t\om-\De\om|\leq c(\A)(|\om|+|\na\om|)\ee for any
$|y'|>20$ and for any $s\in [-20,0]$. Moreover, by (\ref{311}),
\be\la{54}\om(\cdot,0)=0\qquad\mbox{in}\quad\rb^3.\ee By the
backward uniqueness results for the heat operator with variable
lower order terms in a half-space, see \ci{ESS3}, \ci{ESS4}, and
\ci{S6}, and, by (\ref{53}) and (\ref{54}),  we state
\be\la{55}\om(y,s)=0\ee
 for any
$|y'|>20$ and for any $s\in [-20,0]$.

Since our solution is sufficiently smooth in $\rb^3\setminus
\{y'\neq 0\}\times [-10,0]$, one can make use of the unique
continuation through spatial boundaries and conclude that
\be\la{56}\na\wedge u\equiv 0\qquad\mbox{ in}\quad \rb^3\setminus
\{y'\neq 0\}\times [-8,0].\ee
 On the other hand, from (\ref{38}), it follows that
$$\A_0\geq \ess \sup_{-20\leq s\leq 0}\int\limits_{|y'|\leq 40}
\frac {|u(y,s)|^2}{|y'|}dy.$$ So, we observe that, for any $s\in S$,
 \be\la{57}\int\limits_{-\infty}^{+\infty}dy_3\int\limits_{|y'|\leq
40}
 \frac {|u(y,s)|^2}{|y'|}dy'\leq \A_0<+\infty,\ee
where $S\subset [-20,0]$ and $|S|=20$.

Now, we wish to show \be\la{58}\na\wedge u(\cdot,s)\equiv
0\qquad\mbox{ in}\quad \rb^3\ee for any $s\in S$. To this end, we
proceed as follows. Let $\varphi\in C_0^\infty(B')$ be a
non-negative cut-off function being equal to 1 in $B'(1/2)$. Here,
$B'$ and $B'(1/2)$ are two-dimensional balls centered at the
origin with radii 1 and 1/2, respectively. Next, let $\psi$ be an
arbitrary smooth, compactly supported in $\rb^3$, vector-valued
function. Then, by (\ref{56}), for any $s\in [-8,0]$,
$$\int\limits_{\rb^3}u(y,s)\cdot\na\wedge\Big(\psi(y)(1-\varphi(y'/R))
\Big)dy=0=$$$$=\int\limits_{\rb^3}u(y,s)\cdot\na\wedge\psi(y)dy-
\int\limits_{\rb^3}u(y,s)\cdot\na\wedge\Big(\psi(y)\varphi(y'/R)\Big)dy
=J_1(s)+J_2(s).$$ For $J_2$, we have the estimate
$$|J_2(s)|\leq c\Big(1+\frac 1R\Big)\int\limits_{\spt \psi\cap\{|y'|<R\}}
|u(y,s)|dy=$$$$=c\Big(1+\frac 1R\Big)\int\limits_{\spt
\psi\cap\{|y'|<R\}} \frac {|u(y,s)|}{|y'|^\frac 12}|y'|^\frac
12dy\leq $$
$$\leq c\Big(1+\frac 1R\Big)\Big(\int\limits_{\spt
\psi\cap\{|y'|<R\}} \frac {|u(y,s)|^2}{|y'|}dy\Big)^\frac
12\Big(\int\limits_{\spt \psi\cap\{|y'|<R\}} {|y'|}dy\Big)^\frac
12\leq$$$$\leq c(\psi)\Big(1+\frac 1R\Big)R^\frac
32\Big(\int\limits_{-\infty}^{+\infty}dy_3\int\limits_{|y'|<40}\frac
{|u(y,s)|^2}{|y'|}dy'\Big)^\frac 12.$$ By (\ref{57}), the right hand
side of the latter inequality goes to zero as $R\to 0$ for any $s\in
S$. Hence, $J_1(s)=0$ for any $s\in S\cap [-8,0]$, which is but a
weak form of (\ref{58}). By the fact that $u$ is divergence free, we
then show
$$\De u(\cdot,s)=0\qquad \mbox{in}\quad \rb^3$$
for any $s\in S\cap [-8,0]$.

Now, let $B(y_0,R)$ be a ball of radius $R$ with the center at the
point $y_0$. For any $y_0\in \{|y'|\leq 30,\,y_3\in \rb\}$,
$$B(y_0,1)\subset\{|y'|\leq 40,\,y_3\in \rb\}$$ and, since $u$ is
 harmonic,
 $$|u(y_0,s)|\leq
 c\Big(\int\limits_{B(y_0,1)}|u(y,s)|^2dy
 \Big)^\frac 12\leq c\Big(\int\limits_{|y'|\leq 40}|u(y,s)|^2dy
 \Big)^\frac 12\leq $$$$\leq c\sqrt{40\A_0}$$
for any $s\in S\cap [-8,0]$. So, according to (\ref{52}), the
function $u(\cdot,s)$ is bounded in $\rb^3$ for any $s\in S\cap
[-8,0]$. But, by (\ref{51}), in fact, $u(\cdot,s)=0$ in $\rb^3$ for
any $s\in S\cap [-8,0]$. This contradicts with (\ref{39}). Theorem
\ref{1t2} is proved.

G. Seregin\\
Steklov Institute of Mathematics at St.Petersburg, \\
St.Petersburg, Russia
\\
\\
W.Zajaczkowski\\
Institute of Mathematics, Polish Academy of Sciences, Sniadeckich 8,
00-956 Warsaw, Poland

\begin{thebibliography}{99}

\bibitem {CKN}
Caffarelli, L., Kohn, R.-V., Nirenberg, L., Partial regularity of
suitable weak solutions of the Navier-Stokes equations, Comm. Pure
Appl. Math., Vol. XXXV (1982), pp. 771--831.

\bibitem{CL}
Chae D., Lee, J., On the regularity of the axisymmetric solutions
of the Navier-Stokes equations, Math. Z., 239(2002), 645-671.

\bibitem {ESS3}

Escauriaza,L., Seregin, G.,  ~\v Sver\'ak, V., Backward uniqueness
for the heat operator in half space,  Algebra and Analyis,
15(2003), no. 1, 201-214.
\bibitem{ESS4}
Escauriaza,L., Seregin, G.,  ~\v Sver\'ak, V.,.
$L_{3,\infty}$-Solutions to the Navier-Stokes equations and
backward uniqueness, Uspekhi Matematicheskih Nauk, v. 58, 2(350),
pp. 3--44. English translation in Russian Mathematical Surveys,
58(2003)2, pp. 211-250.
\bibitem {Gi} Giga, Y., Solutions for semilinear parabolic
equations in $L^p$ and regularity of weak solutions of the
Navier-Stokes equations, J. of Diff. Equations, 62(1986), pp.
186--212.


\bibitem {L1a}
Ladyzhenskaya, O. A., On uniqueness and smoothness of generalized
solutions to the Navier-Stokes equations, Zapiski Nauchn. Seminar.
POMI, 5(1967), pp. 169--185.


\bibitem{L}
Ladyzhenskaya, O. A., On unique solvability of the
three-dimensional Cauchy problem for the Navier-Stokes equations
under the axial symmetry, Zap. Nauchn. Sem. LOMI 7(1968), 155-177.

\bibitem {LS}
Ladyzhenskaya, O. A., Seregin, G. A., On partial regularity of
suitable weak solutions to the three-dimensional Navier-Stokes
equations, J.  math. fluid mech.,  1(1999), pp. 356-387.

\bibitem{LMNP}
Leonardi, S., Malek, Necas, J., \& Pokorny, M., On axially
simmetric flows in $\mathbb R^3$, ZAA, 18(1999),  639-649.

\bibitem {Li}
Lin, F.-H., A new proof of the Caffarelly-Kohn-Nirenberg theorem,
Comm. Pure Appl. Math., 51(1998), no.3, pp. 241--257.

\bibitem{NP} Neustupa, J., Pokorny, M., Axisymmetric flow of
Navier-Stokes fluid in the whole space with non-zero angular
velocity compnents, Math. Bohemica, 126(2001), 469-481.
\bibitem{Po}
Pokorny, M., A regularity criterion for the angular velocity
component in the case of axisymmetric Navier-Stokes equations,
2001.
\bibitem {Pr}
Prodi, G., Un teorema di unicit\`a per el equazioni di
Navier-Stokes, Ann. Mat. Pura Appl., 48(1959), pp. 173--182.

\bibitem {S2}
Seregin, G. A. On the number of singular points of weak solutions
to the Navier-Stokes equations, Comm. Pure Appl. Math., 54(2001),
issue 8, pp. 1019-1028.
\bibitem {S6}
 Seregin, G.A., On smoothness of $L_{3,\infty}$-solutions
 to the Navier-Stokes equations up to boundary,
 Mathematische Annalen,  332(2005), pp. 219-238.


\bibitem{S7} Seregin, G., New version of the
Ladyzhenskaya-Prodi-Serrin condition, Algebra i Analiz, 18(2006),
124-143.

\bibitem{S8} Seregin, G., Local regularity theory of the
Navier-Stokes equations, to appear in Handbook of Mathematical
Fluid Mechanics, vol. 4.

\bibitem {S9} Seregin, G., Estimates of suitable weak solutions
to the Navier-Stokes equations in critical Morrey spaces, Zapiski
Nauch. Seminar POMI, 2006.

\bibitem{Se}
Serrin, J.,  On the interior regularity of weak solutions of the
Navier-Stokes equations, Arch. Ration. Mech. Anal., 9(1962), pp.
187--195.

\bibitem{Str0}

Struwe, M., On partial regularity results for the Navier-Stokes
equations, Comm. Pure Appl. Math., 41 (1988), no. 4, 437--458.


\bibitem{UY}
Ukhovskij, M. R., Yudovich, V. L., Axially symmetric motions of
ideal and viscous fluids filling all space, Prikl. Mat. Mech. 32
(1968), 59-69.

\bibitem {WZ}
Wiegner, M., Zajaczkowski, W. M., On stability of axially
symmetric solutions to Navier-Stokes equations in a cylindrical
domain and with boundary slip conditions, Banach Center Publ.,
70(2005), 251-278.

\bibitem {Z1}
Zajaczkowski, W. M., Global special regular solutions to the
Navier-Stokes equations in axially symmetric domains under
boundary slip conditions, Diss. Math., 400(2005).


\bibitem{Z2}
Zajaczkowski, W. M., Global special regular to the Navier-Stokes
equations in a cylindrical domain under boundary slip conditions,
Gakuto Series in Mathematics, vol. 21, 2004.








\end{thebibliography}
\end{document}